\documentclass[a4paper,12pt]{article}
\usepackage{amssymb,amsmath,amsthm,hyperref}
\usepackage[mathscr]{eucal}

\newtheorem{prop}{Proposition}[section]
\newtheorem{cor}[prop]{Corollary}
\newtheorem{theorem}[prop]{Theorem}
\newtheorem{lemma}[prop]{Lemma}

\theoremstyle{definition}
\newtheorem{defn}[prop]{Definition}

\newtheorem{examples}[prop]{Examples}

\newenvironment{tightlist}{\begin{list}{$\bullet$}
{\setlength{\itemsep}{0ex plus0.2ex}
\setlength{\topsep}{0ex plus0.2ex}}}
{\end{list}}

\DeclareMathOperator{\Mod}{Mod}  
\DeclareMathOperator{\qftp}{qftp}  

\DeclareMathOperator{\Aut}{Aut}  
\DeclareMathOperator{\cl}{cl}   

\DeclareMathOperator{\im}{Im}   
\DeclareMathOperator{\preim}{preim} 

\newcommand{\card}[1]{\left| #1 \right|}   

\newcommand{\entails}{\ensuremath{\vdash}}   

\newcommand{\restrict}[1]{\ensuremath{\!\!\upharpoonright_{#1}}}

\newcommand{\N}{\ensuremath{\mathbb{N}}}
\newcommand{\Z}{\ensuremath{\mathbb{Z}}}

\newcommand{\ON}{\ensuremath{\mathrm{ON}}}  

\newcommand{\Cat}{\ensuremath{\mathcal{C}}} 
\renewcommand{\L}{L} 
\newcommand{\Loo}{\ensuremath{\L_{\omega_1,\omega}}} 
\newcommand{\Looq}{\ensuremath{\Loo(Q)}}

\renewcommand{\phi}{\varphi}

\renewcommand{\le}{\ensuremath{\leqslant}}
\renewcommand{\ge}{\ensuremath{\geqslant}}
\newcommand{\tuple}[1]{\ensuremath{\langle #1 \rangle}}
\newcommand{\class}[2]{\ensuremath{\left\{ #1 \,\left|\, #2 \right.\right\}}}

\newcommand{\iso}{\cong}

\newcommand{\into}{\hookrightarrow}

\newcommand{\pemb}{\ensuremath{\rightharpoonup}} 

\newcommand{\subs}{\subseteq} 
\newcommand{\prsubs}{\subsetneq}  

\DeclareMathOperator{\fin}{fin}  
\newcommand{\subsetfin}{\ensuremath{\subseteq_{\fin}}} 
\newcommand{\finsub}{\subsetfin}

\newcommand{\minus}{\ensuremath{\smallsetminus}}
\newcommand{\powerset}{\ensuremath{\mathcal{P}}} 

\newcommand{\strong}{\ensuremath{\lhd}} 
\newcommand{\nstrong}{\ensuremath{\not\kern-4pt\lhd\;}} 

\newcommand{\ra}[3]{\ensuremath{#1 \stackrel{#2}{\longrightarrow} #3}}

\newcommand{\leteq}{\mathrel{\mathop:}=}

\DeclareMathOperator{\Indep}{Indep}

\renewcommand{\strong}{\preccurlyeq}

\newcommand{\D}{\ensuremath{\mathcal{D}}}

\newcommand{\qmeclass}{quasiminimal excellent class}

\title{On Quasiminimal Excellent Classes}
\author{Jonathan Kirby\thanks{University of Oxford and University of
    Illinois at Chicago, supported by the EPSRC}}
\date{Version 2.3, 25${}^{\mathrm{th}}$ October 2007}

\begin{document}

\maketitle

\begin{abstract}
  A careful exposition of Zilber's quasiminimal excellent classes and
  their categoricity is given, leading to two new results: the
  $L_{\omega_1,\omega}(Q)$-definability assumption may be
  dropped, and each class is determined by its model of dimension
  $\aleph_0$.
\end{abstract}

Boris Zilber developed quasiminimal excellent classes in
\cite{Zilber05qmec} in order to prove that his conjectural description
of complex exponentiation was uncountably categorical, that is, it has
exactly one model of each uncountable cardinality. This article gives
a simplified and careful exposition of \qmeclass es, and of the
categoricity proof. This more careful exposition has led to two new
results. We say that a \qmeclass\ is \emph{degenerate} iff either it has only
finite dimensional models, or it is a proper subclass of another
\qmeclass. These are essentially uninteresting cases (at least in the
context for which \qmeclass es were invented).  In
\cite{Zilber05qmec}, the proof of the existence of arbitrarily large
models depended on the class being definable by an $\Looq$-sentence of
a specific form (and on having an infinite-dimensional model). The
question of whether this could be generalized to any $\Looq$-sentence
was posed. Here we show that every nondegenerate \qmeclass\ is
definable by an $\Looq$-sentence of the specific form given, hence is
uncountably categorical. Furthermore, any \qmeclass\ with an
infinite-dimensional model extends uniquely to a nondegenerate class,
and hence a \qmeclass\ may be produced in a ``bootstrap'' fashion from
its unique model of dimension $\aleph_0$.

This article sprang from many lively and productive discussions I had
with John Baldwin. The account of \qmeclass es in
\cite{Baldwin_monograph} was rewritten in conjunction with this
article, and incorporates many of the ideas given here. Some of the
differences between this account of the categoricity proof and Boris
Zilber's are due to John Baldwin, in particular the construction of
the isomorphism in theorem~\ref{combinatorial principle} as the union
of the maps $f_X$, and the introduction of the Shelah-style statement
of excellence in lemma~\ref{excellence lemma}.

\section{The Definition}

\begin{defn}
 
  A \emph{quasiminimal excellent class} consists of the following data,
  satisfying axioms I, II, and III. \\

\noindent \textbf{Data:}
  \begin{itemize}
    \item For a given first-order language $\L$, a class $\Cat$ of
      $\L$-structures.
    \item For each $H \in \Cat$, a function $\ra{\powerset
        H}{\cl_H}{\powerset H}$.
  \end{itemize}

\noindent \textbf{Axioms:}
  \begin{description}
    \item[I: Pregeometry] 
      \begin{tightlist}
        \item[]
        \item[I.1] For each $H \in \Cat$, $\cl_H$ is a pregeometry on
          $H$, satisfying the countable closure property (CCP): the
          closure of any finite set is countable.
	\item[I.2] If $H \in \Cat$ and $X \subs H$, then $\cl_H(X) \in
	  \Cat$.
	\item[I.3] If $H \in \Cat$, $X \subs H$, $y \in \cl_H(X)$, and
	  $f: H \pemb H'$ is a partial embedding with $X \cup \{y\}
	  \subs \preim(f)$, then $f(y)\in \cl_{H'}(f(X))$.
      \end{tightlist}

    \item[II: $\aleph_0$-homogeneity over countable models]

      Let $H,H'\in \Cat$, let $G \subs H$ and $G' \subs H'$ be
	  countable closed subsets or the empty set,
	  and let $g: G\to G'$ be an isomorphism.
      \begin{tightlist}
        \item[]
        \item[II.1] If $x \in H$ and $x' \in H'$, independent from $G$
	  and $G'$ respectively, then $g \cup \{(x,x')\}$ is a partial
	  embedding.

	\item[II.2] If $g \cup f: H \pemb H'$ is a partial embedding,
	  $f$ has finite preimage $X$, and $y \in \cl_H(X\cup G)$,
	  then there is $y' \in H'$ such that $g \cup f \cup
	  \{(y,y')\}$ is a partial embedding.
      \end{tightlist}
  \end{description}

  \noindent A \emph{crown} in $H$ is a subset $C \subs H$ such that
    there is an independent subset $B$ of $H$ and finitely many
    subsets $B_1,\ldots,B_n$ of $B$ such that $C = \bigcup_{i=1}^n
    \cl_H B_i$.

  \begin{description}
    \item[III: Quasiminimal excellence]\ \\ 
      Let $H,H' \in \Cat$, let $C$ be a countable crown in $H$, and
      let $g: H \pemb H'$ be a closed partial embedding defined on
      $C$. For any finite subset $X$ of $\cl_H(C)$, there is a finite
      subset $C_0$ of $C$ such that if $f: H \pemb H'$ has preimage
      $X$ and $f \cup g\restrict{C_0}$ is a partial embedding then $f
      \cup g$ is also a partial embedding. We say that the
      quantifier-free type of $X$ over $C$ is \emph{determined over
      $C_0$}.
  \end{description}
\end{defn}

\begin{defn}
  We consider the class $\Cat$ as a category by taking the
  \emph{closed} embeddings, that is, those $\L$-embeddings $H \into
  H'$ such that the image of $H$ is closed in $H'$. Write $H \strong
  H'$ if the inclusion $H \subs H'$ is a closed embedding.

  A partial embedding $f: H \pemb H'$ is \emph{closed} iff for every
  closed set $X$ in the preimage of $f$, the image $f(X)$ is closed in
  $H'$.
\end{defn}

The notion of a closed $\L$-embedding is the right one, because it
also preserves the closure operator.
\begin{lemma}\label{embeddings preserve closure}
  If $H, H' \in \Cat$ with $H \subs H'$ and $X$ is any subset of $H$
  then $\cl_H(X) = \cl_{H'}(X) \cap H$. Furthermore, if $H \strong H'$
  then $\cl_H(X) = \cl_{H'}(X)$. In particular, a closed embedding is
  a closed partial embedding (with respect to the definitions above).
\end{lemma}
\begin{proof}
  The first statement is axiom I.3 applied to the inclusion map $H
  \into H'$. For the second statement, since $X \subs H$ we have
  $\cl_{H'}(X) \subs \cl_{H'}(H) = H$, and so $\cl_{H'}(X) \cap H =
  \cl_{H'}(X)$.
\end{proof}

The definition of a \qmeclass\ given here differs from Zilber's in the
following ways. The pregeometry is included as data rather than its
existence being postulated as an axiom. This avoids ambiguity where
there may be more than one pregeometry to choose from. The exchange
and CCP axioms are included in the definition. Early versions of
Zilber's axiomatization omitted exchange, but he later realized it was
necessary. In the original definition, only the models with CCP were
of interest (indeed only they were quasiminimal), but for technical
reasons CCP was not included as an axiom. The technical reasons are
avoided by considering a \qmeclass\ as a category with closed
embeddings (which are also introduced here). This is very natural, as
it makes the class into an abstract elementary class (provided that it
has unions of chains).  Axiom I.3 was listed among the homogeneity
axioms rather than the pregeometry axioms in \cite{Zilber05qmec}, and
is here made precise.  Axiom II is now $\aleph_0$-homogeneity over
\emph{countable, closed} submodels, not all submodels. II.1 is
weakened to consider only singletons, not independent tuples of
arbitrary length. The excellence axiom III is stated only for
countable models. These weakenings of the axioms to countable models
and the corresponding strengthening of the categoricity theorem are
the main reason why the new results can be proved.

The terminology of a type being \emph{determined} over a set is
preferred to Zilber's original \emph{defined} over a set, which
conflicts with another, different usage. I have introduced the
terminology \emph{crown} which is shorter than \emph{union of an
  independent $n$-system}, or \emph{independent $n$-cube}, does not
contain the awkward parameter $n$, and which I think is a better
description of the concept. Zilber has also used \emph{special subset}.

The axioms all refer to partial embeddings, and imply that the
language $\L$ is rich enough to have a form of quantifier elimination
(see proposition~\ref{types}). This is a minor convenience, but is not
in any way necessary. In examples, this quantifier elimination must
usually be obtained by first expanding the language. It would perhaps
be better to define a quasiminimal excellent class to be any class
$\Cat'$ of $\L'$-structures for which there is an expansion by
definitions to a language $\L$, such that the resulting class $\Cat$
satisfies the given axioms. Alternatively, one could work throughout
with closed partial $\Loo$-maps, and adjust the axioms accordingly.

In section~\ref{existence section}, I introduce additional axioms IV
on unions of chains and the existence of an infinite-dimensional
model.  The subsequent analysis shows that adding these axioms to the
definition of a \qmeclass\ rules out exactly the degenerate cases,
hence it would be natural and would do no harm, but that convention is
not adopted in this paper.

The reader may like to have some examples in mind. The third of the
examples below is the simplest in which not all submodels are closed.
\begin{examples}
The following are \qmeclass es.
  \begin{itemize}
    \item Any strongly minimal theory in a countable language, with
      algebraic closure. (These are elementary classes.)
    \item $\mathrm{ACF}_0$, with a predicate $Z$ and the axiom
      $Z(x) \longleftrightarrow \bigvee_{n\in \Z} x = n$, with
      algebraic closure (an $\Loo$-class).
    \item The theory consisting just of one equivalence relation, all
      of whose blocks have size $\aleph_0$, the closure of a subset
      $X$ being the union of the blocks meeting $X$ (an $\Looq$-class).
    \item The motivating example: Zilber's exponential field, and the
      related examples of ``covers of the multiplicative group'' and
      ``raising to powers''.
  \end{itemize}
\end{examples}

\section{Uniqueness of models up to dimension $\aleph_1$}

\begin{theorem}\label{combinatorial principle up to aleph1}
  Let $\Cat$ be a class satisfying axioms I and II, and let $H,H' \in
  \Cat$ with $\dim H \le \aleph_1$. Let $G \subs H$ be empty or closed
  and countable, and let $f_0:G \to H'$ be a closed embedding (or the
  empty map if $G$ is empty). Let $B$ be a basis of $H$ over $G$ and
  suppose $\psi_B: H \pemb H'$ is an injective partial map with
  preimage $B$ and image an independent set over $\im f_0$. Then $\psi
  \leteq f_0 \cup \psi_B$ extends to a closed map $\hat{\psi}: H \to
  H'$.

  In particular, if $\im{\psi}$ spans $H'$ then $\hat{\psi}$ is an
  isomorphism.
\end{theorem}
\begin{proof}
  Well-order $B$ as $(b_\lambda)_{\lambda<\mu}$ for some ordinal $\mu
  \le \omega_1$. For each ordinal $\nu \le \mu$, let $G_\nu = \cl_H(G
  \cup \class{b_\lambda}{\lambda<\nu})$. In particular, $G_0 =
  G$. Inductively, we construct closed maps $f_\nu : G_\nu \to H'$ such
  that 
  \begin{tightlist}
    \item $\psi\restrict{G_\nu} \subs f_\nu$, and
    \item If $\nu_1 \le \nu_2$ then $f_{\nu_1} \subs f_{\nu_2}$.
  \end{tightlist}
  At limit ordinals, take unions. For a successor $\nu = \lambda+1$,
  we construct $f_\nu$ as $\bigcup_{n\in\N} h_n$ where the $h_n$ are
  partial embeddings constructed inductively by the back and forth
  method, such that $h_0 = \{(b_\lambda,\psi(b_\lambda))\}$, and for
  each $n$,
  \begin{tightlist}
    \item $\preim(h_n)$ is finite,
    \item $h_n \subs h_{n+1}$, and
    \item $h_n \cup f_\lambda$ is a partial embedding.
  \end{tightlist}

  Both $G_\nu$ and $G_\nu' \leteq
  \cl_{H'}(f_\lambda(G_\lambda)\psi(b_\lambda))$ are countable, so
  list their elements in chains of length $\omega$.

  For odd $n$, let $a$ be the least element of $G_\nu \minus \preim
  h_{n-1}$. Using the $\aleph_0$-homogeneity over the countable model
  $G_\lambda$, there is $b \in G_\nu'$ such that, taking $h_n \leteq
  h_{n-1} \cup \{(a,b)\}$, $h_n \cup f_\lambda$ is a partial
  embedding.

  For even $n$, repeat the process going back rather than forth. Then
  $f_\nu$ is a embedding defined on all of $G_\nu$, with image all of
  $G_\nu'$, hence it is a closed embedding. It extends $f_\lambda$ and
  $\psi\restrict{G_\nu}$ by construction.

  Take $\hat{\psi} = \bigcup_{\lambda<\mu}f_\lambda$. This is a closed
  embedding as required.
\end{proof}

\begin{cor}[Uniqueness of models up to $\aleph_1$]\label{uniqueness up
    to aleph1} 
  Models of a quasiminimal excellent class of dimension less than or
  equal to $\aleph_1$ are determined up to isomorphism by their
  dimension. There is at most one model of cardinality $\aleph_1$.
  \qed
\end{cor}

A small modification of the proof explains why axiom II is called
$\aleph_0$-homogeneity over countable models.
\begin{prop}\label{aleph_0 homogeneity}
  Let $\Cat$ be a class satisfying axioms I and II, let $H \in \Cat$,
  let $G \subs H$ be empty or countable and closed, and let
  $\bar{x},\bar{y}$ be $n$-tuples from $H$ with the same
  quantifier-free type over $G$. Then there is an automorphism $\sigma
  \in \Aut(H/G)$ such that $\sigma(\bar{x}) = \bar{y}$.
\end{prop}
\begin{proof}
  Write the tuple $\bar{x}$ as $b_1,\ldots,b_n,a_1,\ldots,a_m$ where
  the $b_i$ are independent over $G$ and each $a_i$ lies in $\cl_H(G
  b_1,\ldots,b_n)$, and correspondingly write $\bar{y}$ as
  $b_1',\ldots,b_n',a_1',\ldots,a_m'$. Extend $\bar{b}$ and $\bar{b}'$
  to bases of $H$ over $G$, and let $\psi_B$ be a bijection between
  the bases which sends $b_i$ to $b_i'$ for $i=1,\ldots,n$. Let $f_0$
  be the identity map on $G$. Then follow the proof of
  theorem~\ref{combinatorial principle up to aleph1}, except that in
  the construction of $f_\nu$ for $\nu \le n$, start by sending $a_i$
  to $a_i'$ for each $i$ such that $a_i \in \cl_H(Gb_1,\ldots,b_\nu)$.
\end{proof}

Note that the construction in the proof of theorem~\ref{combinatorial
  principle up to aleph1} cannot be carried out for $H$ of dimension
greater than $\aleph_1$, because we only have $\aleph_0$-homogeneity
over models, not $\aleph_1$-homogeneity over models. (In addition, we
have only assumed homogeneity over countable models, but we show later
that is not an issue, whereas most examples are actually not
$\aleph_1$-homogeneous over models.)

\section{Uniqueness of large models}

\begin{lemma}\label{compatibility lemma}
  Let $M$ be a set, let $\cl$ be a pregeometry on $M$ and let $B$ be
  an independent subset of $M$. Let $X, Y \subs B$. Then $\cl(X) \cap
  \cl(Y) = \cl(X\cap Y)$.
\end{lemma}
\begin{proof}
  By monotonicity, $\cl(X \cap Y) \subs \cl(X)$ and $\cl(X \cap Y)
  \subs \cl(Y)$, so $\cl(X \cap Y) \subs \cl(X)\cap \cl(Y)$.

  Now suppose $z \in \cl(X) \cap \cl(Y) \minus \cl(X \cap Y)$. By
  finite character, there are $X_0 \finsub X$ and $Y_0 \finsub Y$ such
  that $z \in \cl(X_0) \cap \cl(Y_0) \minus \cl(X \cap Y)$. By
  exchange, there is $x\in X_0 \minus Y$ such that $\cl(X_0) = \cl(X_0
  z - x)$. Similarly, there is $y \in Y_0 \minus X$ such that
  $\cl(Y_0) = \cl(Y_0 z - y)$. Hence
  \[C \leteq \cl(X_0 Y_0) = \cl (X_0 Y_0 z -x -y)\]
  but $X_0 \cup Y_0$ is an independent set and so $\dim C = \card{X_0
  \cup Y_0}$, and yet $C$ is generated by a set of size $\card{X_0
  \cup Y_0} - 1$, which is a contradiction. Hence no such $z$ exists,
  and $\cl(X) \cap \cl(Y) = \cl(X\cap Y)$ as required.
\end{proof}

We translate Zilber's excellence criterion (which can be proved
directly in examples) to a Shelah-style criterion.
\begin{lemma}[Excellence -- Shelah style]\label{excellence lemma}
  Let $H,H' \in \Cat$ and let $C \subs H$ be a countable crown. Then
  every closed partial embedding $f: H \pemb H'$ with preimage $C$
  extends to a closed embedding $\hat{f}: \cl_H(C) \into H'$.
\end{lemma}
\begin{proof}
  Let $\bar{C} = \cl_H(C)$ and let $\bar{C}' = \cl_{H'}(f(C))$. They
  are both countable, so choose an ordering of each of length
  $\omega$.

  Inductively we construction partial embeddings $f_n : H \pemb H'$
  for $n \in \N$ such that for each $n$:
  \begin{tightlist}
    \item $\preim(f_n)$ is finite,
    \item $f_n \subs f_{n+1}$, and
    \item $f_n \cup f$ is a partial embedding.
  \end{tightlist}
  Take $f_0 = \emptyset$. We construct the $f_n$ for $n > 0$ via the
  back and forth method, going forth for odd $n$ and back for even
  $n$.

  For odd $n$, let $a$ be the least element of $\bar{C} \minus \preim
  f_{n-1}$. Then the set $\preim f_{n-1} \cup\{a\}$ is a finite subset
  of $\bar{C}$, so by axiom III (excellence) and the finite character
  of the pregeometry there is a finite subset $C_0$ of $C$ such that
  the quantifier-free type of $\preim f_{n-1} \cup\{a\}$ over $C$ is
  determined over $C_0$ and $a \in \cl_H(C_0)$. Let $g = f_{n-1} \cup
  f\restrict{C_0}$. By induction, $f_{n-1} \cup f$ is a partial
  embedding, so $g$ is a partial embedding. By axiom II.2
  ($\aleph_0$-homogeneity), there is $b \in H'$ such that $f_n \leteq
  g \cup \{(a,b)\}$ is a partial embedding. Since the type of $\preim
  f_n$ over $C$ is determined over $C_0$, $f_n \cup f$ is a partial
  embedding, as required.

  For even $n$, note that $f(C)$ is a crown in $H'$ because $f$ is a
  closed partial embedding. Also note that the inverse of a partial
  embedding is a partial embedding. Hence we can perform the same
  process as for odd steps, reversing the roles of $H$ and $H'$, to
  find $f_n$ whose image contains the least element of $\bar{C}'$ not
  in the image of $f_{n-1}$.

  Let $\hat{f} = \bigcup_{n \in \N} f_n$. Then $\hat{f}$ is an
  embedding extending $f$, defined on all of $\bar{C}$, whose image is
  all of $\bar{C}'$. Hence $\hat{f}$ is a closed embedding.
\end{proof}

\begin{theorem}\label{combinatorial principle} 
  Let $\Cat$ be a quasiminimal excellent class, and let $H,H' \in
  \Cat$. Let $G\subs H$ be empty or closed and let $f_0: G \to H'$ be
  a closed embedding. Let $B$ be a basis of $H$ over $G$ and suppose
  $\psi_B: H \pemb H'$ is an injective partial map with preimage $B$
  and image an independent set over $\im f_0$. Then $\psi \leteq f_0
  \cup \psi_B$ extends to a closed embedding $\hat{\psi}: H \to H'$.

  In particular, if $\im{\psi}$ spans $H'$ then $\hat{\psi}$ is an
  isomorphism.
\end{theorem}

\begin{proof}
  If $\dim H$ is finite then we are done by theorem~\ref{combinatorial
  principle up to aleph1}, so we assume $\dim H$ is infinite. If $G$
  is empty or finite dimensional then we may extend $G$ to a closed
  subset of dimension $\aleph_0$, using theorem~\ref{combinatorial
  principle up to aleph1}. We first prove the theorem assuming that
  $\dim G = \aleph_0$. The case where $\dim G > \aleph_0$ will be
  discussed afterwards.

  For each finite subset $X \finsub B$, we will construct a closed
  embedding $f_X: \cl_H(GX) \to H'$ such that
  \begin{tightlist}
    \item whenever $Y \subs X$, we have $f_Y \subs f_X$, and
    \item $f_X \restrict{X} = \psi\restrict{X}$. 
  \end{tightlist}
  We construct the $f_X$ by well-founded induction on the partial order
  of finite subsets of $B$. Take $f_\emptyset \leteq f_0$.

  Suppose that $X\finsub B$, $X\neq \emptyset$, and we already have
  $f_Y$ for all proper subsets $Y$ of $X$. Let $g_X = \bigcup_{Y
    \prsubs X} f_Y$. Using lemma~\ref{compatibility lemma}, we see
  that if $x \in \preim f_{Y_1} \cap \preim f_{Y_2}$ for $Y_1,Y_2
  \prsubs X$ then $x \in \preim f_{Y_1 \cap Y_2}$, and so by
  hypothesis $f_{Y_1}(x) = f_{Y_2}(x)$. Hence $g_X$ is a well-defined
  partial function. We must show that $g_X$ is a partial embedding.

  Say $X = \{x_1,\ldots, x_n\}$, and for $i=1,\ldots,n$ let $Y_i = X
  \minus \{x_i\}$. Let $C_k = \bigcup_{i=1}^k \cl_H (GY_i)$ and let
  $h_k = \bigcup_{i=1}^k f_{Y_i}$. So $g_X = h_n$.

  We prove by induction on $k$ that $h_k$ is a partial embedding. The
  $k=1$ case is immediate. For the induction step, take tuples $a \in
  C_{k-1}$ and $b \in \cl_H(GY_k)$. We construct an automorphism of a
  submodel of $H'$ which has the effect of moving the parameters $a$
  inside $\cl_H(GY_k)$. Let $A \finsub G$ such that $a,b \in
  \cl_H(AX)$, and let $G_0 = \cl_H(A)$. Let $z \in G \minus G_0$, let
  $H_0 = \cl_H(G_0Xz)$, and let $H_0' = \cl_{H'}(\psi(G_0Xz))$. By
  theorem~\ref{combinatorial principle up to aleph1}, there is an
  automorphism $\sigma$ of $H_0$, fixing $\cl_{H}(G_0Y_k)$ and swapping
  $x_k$ and $z$.

  The idea is to compare $h_k$ on $H_0$ with the composite embedding
  $\tau = \sigma'^{-1}f_{Y_k}\sigma$, where $\sigma'$ is the
  automorphism of $H_0'$ which ``corresponds'' to $\sigma$. To have a
  notion of what ``corresponds'' means, we must choose a suitable
  isomorphism between $H_0$ and $H_0'$. Fortunately, the condition of
  being ``suitable'' is weaker than the compatibility condition we are
  trying to prove!

  Let $e$ be a closed embedding $e : \cl_H(GX) \to H'$ extending
  $h_{k-1} \cup \psi\restrict{X}$. For $k=2$, such an $e$ exists
  because $X = Y_1 \cup \{x_1\}$, so $h_1 \cup \{(x_1,\psi(x_1)\}$
  extends to some $e$ by theorem~\ref{combinatorial principle up to
  aleph1}. For $k \ge 3$, $C_{k-1}$ is a crown whose closure is
  $\cl_H(GX)$, so $h_k$ extends to some $e$ by lemma~\ref{excellence
  lemma}. Let $\sigma' = e\sigma e^{-1}$, and let $\tau = \sigma'^{-1}
  f_{Y_k}\sigma = e\sigma^{-1} e^{-1}f_{Y_k} \sigma$.

  Write $Y_{ik}$ for $Y_i \cap Y_k$. The tuple $a \in
  \bigcup_{i=1}^{k-1}\cl_{H}(G_0Y_i)$, so 
  \[\sigma(a) \in \bigcup_{i=1}^{k-1}\cl_{H}(G_0Y_{ik}z)
  \subs \bigcup_{i=1}^{k-1}\cl_{H}(GY_{ik}) \subs \preim h_{k-1}.\]
  By hypothesis, $f_{Y_k}$ agrees with $f_{Y_i}$ on
  $\cl_{H}(GY_{ik})$, hence $f_{Y_k}$ agrees with $h_{k-1}$ on
  $\bigcup_{i=1}^{k-1}\cl_{H}(GY_{ik})$. Also $e$ and $h_k$ both
  extend $h_{k-1}$, so 
  \[\tau(a) = e \sigma^{-1}e^{-1} f_{Y_k} \sigma(a) 
     = h_{k-1} \sigma^{-1} h_{k-1}^{-1} h_{k-1}\sigma(a) 
     = h_{k-1}(a) = h_k(a).\]

  The tuple $b \in \cl_H(G_0Y_k)$, so it is fixed by $\sigma$. The
  embeddings $f_{Y_k}$ and $e$ preserve the closure, so
  $e^{-1}f_{Y_k}(b)$ is fixed by $\sigma^{-1}$. So 
  \[\tau(b) = e \sigma^{-1} e^{-1} f_{Y_k} \sigma(b) 
     = e e^{-1}f_{Y_k}(b) 
     =  f_{Y_k}(b) = h_k(b).\]
 
  Thus for any quantifier-free formula $R$, 
  \[H \models R(a,b) \iff H' \models R(\tau(a),\tau(b)) \iff H'\models
  R(h_k(a),h_k(b)).\] This holds for any tuples $a,b$ (for a suitable
  choice of $\tau$) and so $h_k$ is a partial embedding. In
  particular, $g_X = h_n$ is a partial embedding. It is a union of
  finitely many closed partial embeddings, hence is a closed partial
  embedding.

  By lemma~\ref{excellence lemma}, $g_X$ extends to a closed embedding
  $f_X : \cl_H(GX) \to H'$. Thus, by induction, we have compatible
  embeddings $f_X$ for every $X \finsub B$. Let $\hat{\psi} =
  \bigcup_{X \finsub B} f_X$, a closed embedding. By the finite
  character of the closure, $H = \bigcup_{X \finsub B}
  \cl_H(GX)$. Hence $\hat{\psi}$ is a total map on $H$. That completes
  the proof in the case where $\dim G \le \aleph_0$.

  If $\dim G > \aleph_0$, let $G'$ be a closed subset of $G$ of
  dimension $\aleph_0$, and let $B_1$ be a basis of $G$ over
  $G'$. Inductively construct embeddings $f_X$ as before, with $G'$ in
  place of $G$ and $B \cup B_1$ in place of $B$, except that for $X
  \finsub B_1$ take $f_X = f_0\restrict{\cl_H(G'X)}$. Then the map
  $\hat{\psi}$ obtained will extend $\psi$ as required.
\end{proof}

\begin{cor}[Uniqueness of models]\label{uniqueness}
  Models of a quasiminimal excellent class are determined up to
  isomorphism by their dimension. There is at most one model of any
  uncountable cardinality. \qed
\end{cor}

A small modification of the proof of theorem~\ref{combinatorial
  principle}, similar to the case of proposition~\ref{aleph_0
  homogeneity}, gives the following full homogeneity and quantifier
elimination statements.
\begin{prop}\label{types}
Let $\Cat$ be a \qmeclass.
\begin{tightlist}
\item Any model in $\Cat$ is $\aleph_0$-homogeneous over the empty set
  and over every closed submodel.
\item The Galois types of finite tuples over the empty set and over
  closed submodels are equal to the quantifier-free $\L$-types. \qed
\end{tightlist}
\end{prop}

\section{Existence of models}\label{existence section}

The axioms for a \qmeclass\ do not imply that any models exist at all.
If there is a model of dimension $\kappa$, then for any cardinal
$\lambda < \kappa$ there is a model of dimension $\lambda$, by axiom
I.2. Hence the models of a quasiminimal excellent class are indexed by
some initial segment of the class of cardinals. There is nothing in
the axioms to say that we cannot have a proper initial segment --
indeed any initial segment of any \qmeclass\ is also a \qmeclass.
There are also \qmeclass es with only models of finite dimension.  To
exclude these degenerate cases, we consider additional axioms. We will
see that they exclude only these degenerate cases.
\begin{description}
  \item[IV: Chains]
    \begin{tightlist}
      \item[]
      \item[IV.1] The category $\Cat$ has unions of chains of all
        ordinal lengths. That is, suppose $(H_\mu)_{\mu<\lambda}$ is
        an ordinal-indexed chain of models of $\Cat$ with closed
        embeddings. Let $H$ be the union of the chain (as an
        $\L$-structure), and for $X \subs H$, define $cl_H(X) =
        \bigcup_{\mu < \lambda} \cl_{H_\mu}(X \cap H_\mu)$. Then
        $\tuple{H,\cl_H} \in \Cat$.
      \item[IV.2] $\Cat$ contains an infinite dimensional model.
    \end{tightlist}
\end{description}

It is easy to show:
\begin{prop}
  A quasiminimal excellent class which satisfies axiom IV.1 (unions of
  chains) is an abstract elementary class with L\"owenheim cardinal at
  most $\aleph_0$.\qed
\end{prop}

\begin{theorem}[Existence of models]\label{existence}
  Let $\Cat$ be a quasiminimal excellent class satisfying axiom IV.
  Then for every cardinal $\kappa$ there is a unique $H \in \Cat$ of
  dimension $\kappa$. In particular, $\Cat$ is uncountably
  categorical. Conversely, any uncountably categorical \qmeclass\
  satisfies axiom IV.
\end{theorem}
\begin{proof}
  We have already proved uniqueness, and it remains to prove
  existence. By induction on ordinals $\lambda$ we construct the
  initial segment of length $\lambda$ of a chain $(M_\mu)_{\mu \in
    \ON}$ in $\Cat$ where each $M_\mu$ has a chosen basis indexed by
  $\mu$, and the inclusion maps $M_\mu \into M_\nu$ of the chain
  extend the inclusion maps $\mu \into \nu$ of the ordinals
  (identified with the chosen bases). Here $\ON$ is the ordered class
  of ordinals.

  By axiom IV.2, there is $M \in \Cat$ with $\dim M = \kappa$, for
  some infinite cardinal $\kappa$. Choose an ordering of a basis of
  $M$ of length the initial ordinal $\alpha$ of size $\kappa$. Taking
  the closures of the initial segments of this basis gives the chain
  $(M_\mu)_{\mu < \alpha}$ in $\Cat$.

  There are two cases for the inductive step. If $\lambda$ is a
  successor $\lambda = \nu^+$ then we already have a model (for
  example $M_\nu$) of dimension $\card{\lambda}$. We choose a new
  ordering of length $\lambda$ for a basis, and, using
  theorem~\ref{combinatorial principle}, we choose a closed embedding
  $M_\nu \into M_\lambda$ extending the inclusion of bases.

  If $\lambda$ is a limit ordinal then by induction we have a chain
  $(M_\mu)_{\mu < \lambda}$ in $\Cat$. By axiom IV.1, the union of
  this chain lies in $\Cat$, and the inclusion maps of the $M_\mu$
  into the limit are closed. It has a basis indexed by $\lambda$, so
  we take it as $M_\lambda$. That gives one direction.

  For the converse, suppose $\Cat$ is uncountably categorical, and
  $(H_\nu)_{\nu < \lambda}$ is a chain in $\Cat$. Let $\kappa$ be a
  cardinal which is an upper bound for the dimensions of the $H_\nu$,
  and for $\card{\lambda}$, and let $\alpha$ be the initial ordinal of
  $\kappa$. Then the model $M_\kappa$ in $\Cat$ of dimension $\kappa$
  can be written as the union of a chain $(M_\mu)_{\mu < \alpha}$, as
  above. The chain $(H_\nu)_{\nu < \lambda}$ is isomorphic to a
  subchain of $(M_\mu)_{\mu < \alpha}$ (possibly with repeats), so its
  union is naturally a closed subset of $M_\kappa$. Hence the union
  lies in $\Cat$. Thus $\Cat$ satisfies IV.1. Axiom IV.2 is immediate.
\end{proof}

\section{Definability}

This section is mainly about the definability of a \qmeclass, but we
first give a result about definable sets within a class. The
motivation for the word \emph{quasiminimal} is that every definable
set is countable or cocountable.
\begin{lemma}[Quasiminimality]\label{qm Galois}
  Let $\Cat$ be a quasiminimal excellent class, and let $H \in \Cat$.
  \begin{tightlist}
  \item If $X \subs H$ and $a, b \in H \minus \cl_H(X)$ then there is
    $\sigma \in \Aut(H/X)$ such that $\sigma(a) = b$.
  \item Every subset of $H$ which is definable with countably many
    parameters in any logic respecting $\L$-automorphisms (for example
    $\L_{\infty,\omega}$ or $\Loo(Q)$) is either countable or
    cocountable.
  \end{tightlist}
\end{lemma}
\begin{proof}
  The first part is immediate from theorem~\ref{combinatorial
    principle}. For the second part, let $S \subs H$ be definable from
  a countable parameter set $A$. Then $H\minus \cl_H(A)$ is a single
  orbit of $\Aut(H/A)$, and is cocountable, and either $S$ does not
  meet $H\minus \cl_H(A)$ or $S$ contains $H\minus \cl_H(A)$.
\end{proof}

There is no assumption that the language $\L$ should be countable, but
in fact we may assume that it is.
\begin{prop}\label{countable language}
  Let $\Cat$ be a quasiminimal excellent class, in a language $\L$.
  Then there is a countable sublanguage $\L^0$ of $\L$ such that the
  class $\Cat^0$ of reducts to $\L^0$ of models in $\Cat$ is also a
  quasiminimal excellent class, and the reduct map $\Cat \to \Cat^0$
  is an isomorphism of categories.
\end{prop}
\begin{proof}
  Let $M$ be a countable model, $G$ a closed submodel of $M$ or
  $\emptyset$, and $a,b$ be $n$-tuples from $M$, and suppose the
  quantifier-free types $\qftp(a/G)$ and $\qftp(b/G)$ are different.
  Then there is a symbol from the signature of $\L$ which witnesses
  this. Up to isomorphism, there are only countably many such tuples
  $\tuple{M,G,a,b}$. Hence there is a countable sublanguage $\L^0$ of
  $\L$ which witnesses all such differences of quantifier-free types.
  By induction using $\aleph_0$-homogeneity over models, partial
  $\L^0$-embeddings and partial $\L$-embeddings of countable models
  coincide.

  The language $\L$ is finitary (that is, every symbol has finite
  arity), so a map between two models is an $\L$-embedding precisely
  when its restriction to every finitely generated substructure is an
  $\L$-embedding. By the countable closure property, an embedding is
  closed precisely when its restriction to every countable submodel is
  closed. Hence partial (and total) closed $\L^0$-embeddings and
  closed $\L$-embeddings of any models coincide. In particular,
  quantifier-free $\L^0$-types coincide with quantifier-free
  $\L$-types. Hence for each relation symbol $R$ of $\L$, there is an
  quantifier-free $\L^0$-formula $\theta_R(\bar{x})$ such that for any
  $M \in \Cat$, $M \models (\forall \bar{x})[R(\bar{x})
  \leftrightarrow \theta_R(\bar{x})]$.  This gives a unique way of
  expanding a model in $\Cat^0$ to a model in $\Cat$, which in turn
  gives an inverse to the reduct map.
\end{proof}

\begin{lemma}[$\Loo$-definability of the pregeometry]
  Let $\Cat$ be a \qmeclass. For each $n \in \N$ there is a
  quantifier-free $\Loo$-formula $\pi_n(x,y_1,\ldots,y_n)$ such that
  for each $H \in \Cat$ and each $a,b_1,\ldots,b_n \in H$, we have
  \[a \in \cl_H(\bar{b}) \mbox{ iff } H \models \pi_n(a,\bar{b}).\]
\end{lemma}
\begin{proof}
  By proposition~\ref{countable language} we may assume the language
  is countable. For each $n \in \N$, let $M_n$ be the model of
  dimension $n$ (if it exists, and the model of maximum dimension if
  it does not). For each $n$, every quantifier-free $n$-type is
  realised in the closure of its realization, hence in $M_{n+1}$,
  hence there are only countably many quantifier-free $n$-types over
  $\emptyset$. The types are given by quantifier-free $\Loo$-formulas,
  say $(\beta_j(\bar{y}))_{j \in \N}$. By the same argument, for each
  $j$ and $\bar{b}$ of type $\beta_j$, there are only countably many
  1-types of elements in the closure of $\bar{b}$, and their types are
  also given by $\Loo$-formulas, say $(\gamma_{ij}(x,\bar{b}))_{i \in
    \N}$. So the formula
  \[\pi_n(x,\bar{y}) \equiv \bigvee_{j\in \N} \left(\beta_j(\bar{y})
    \wedge \bigvee_{i \in \N}\gamma_{ij}(x,\bar{y})\right)\] 
  works for $M_{n+1}$. If $H \in \Cat$ and $a,b_1,\ldots,b_n \in H$,
  then either there is a closed embedding $H \into M_{n+1}$, or there
  is a closed embedding $M_{n+1} \into H$ whose image contains
  $a,b_1,\ldots,b_n$. In each case, $H \models \pi_n(a,\bar{b})$ iff
  $M_{n+1} \models \pi_n(a,\bar{b})$ since $\L$-embeddings preserve
  quantifier-free $\Loo$-formulas, and $a \in \cl_H(\bar{b})$ iff $a
  \in \cl_{M_{n+1}}(\bar{b})$ by lemma~\ref{embeddings preserve
    closure}.  Hence the same formulas $\pi_n$ work for every $H \in
  \Cat$.
\end{proof}

\begin{lemma}\label{Loo-embedding}
  Let $\Cat$ be a \qmeclass, let $f: M \to N$ be a closed embedding in
  $\Cat$ and suppose that $\dim M \ge \aleph_0$. Then $f$ is an
  $\Loo$-embedding.
\end{lemma}
\begin{proof}
  We prove by induction on formulas that for any $\Loo$-formula
  $\theta(\bar{x})$ and any $\bar{a} \in M$, we have $M\models
  \theta(\bar{a}) \iff N \models \theta(\bar{a})$.

  \begin{tightlist}
  \item[] $f$ is an $\L$-embedding, so atomic formulas are preserved.
  \item[] The cases $\theta(\bar{x}) \equiv \bigwedge_{i\in I}
    \theta_i(\bar{x})$ and $\theta(\bar{x}) \equiv \neg
    \theta_0(\bar{x})$ are immediate.
  \item[] If $\theta(\bar{x}) \equiv \exists y \phi(y,\bar{x})$ then
    the left to right case is immediate. Suppose $N \models \exists y
    \phi(y,\bar{a})$. Then for some $b \in N$, $N \models
    \phi(b,\bar{a})$. If $b\in M$ then we are done. Suppose $b \notin
    M$. Since $\cl(\bar{a})$ is finite dimensional, we can choose $c
    \in M \minus \cl(\bar{a})$. The point $b$ is independent from $M$,
    so by lemma~\ref{qm Galois} there is an automorphism of $N$ fixing
    $\cl(\bar{a})$ and swapping $b$ and $c$. So $N \models
    \phi(c,\bar{a})$ and, by the inductive hypothesis, $M \models
    \phi(c,\bar{a})$. Hence $M \models \exists y \phi(y,\bar{a})$.
  \end{tightlist}
  Hence $f$ is an $\Loo$-embedding as required.
\end{proof}
It is easy to extend the proof to show that if $\dim H \ge \aleph_1$
then $f$ is an $\Looq$-embedding.

\begin{theorem}\label{Scott sentence}
  Let $\Cat$ be a \qmeclass\ in a countable language, with an infinite
  dimensional model. For each $n \le \omega$, let $M_n$ be the model
  of dimension $n$, and let $\sigma_n$ be its Scott sentence. Let
  \[\Sigma = \left[\bigvee_{n \le\omega} \sigma_n\right] \wedge
  \left[\bigwedge_{n \in \N}(\forall y_1\ldots,y_n) \neg(Qx)
    \pi_n(x,y_1,\ldots,y_n)\right]\] where $Q$ is the quantifier
  ``there exist uncountably many''. Then $\Mod(\Sigma)$ is an
  uncountably categorical \qmeclass\ containing $\Cat$. Furthermore,
  if $\Cat$ satisfies axiom IV then $\Cat = \Mod(\Sigma)$.
\end{theorem}
\begin{proof}
  We check axioms I---IV for $\Mod(\Sigma)$. The statement that the
  $\pi_n$ define a pregeometry can be axiomatized as an
  $\Loo$-sentence, and it is true in each $M_n$, hence it follows from
  each $\sigma_n$. The countable closure property is explicit in
  $\Sigma$, hence axiom I.1 holds. Axiom I.3 holds because the
  pregeometry is defined by the quantifier-free $\Loo$-formulas
  $\pi_n$. From these axioms we get the notion of dimension for each
  model of $\Sigma$, and we also get the notion of a closed embedding.

  For each uncountable cardinal $\kappa$, let $\D_\kappa = \class{N
    \models \Sigma}{\dim N < \kappa}$. Axioms I.1 and I.3 hold for
  each $\D_\kappa$ as well.

  Axioms II and III are statements about the countable models. Any
  countable model of $\Sigma$ must be one of the $M_n$ for $n \le
  \omega$, since Scott sentences are $\aleph_0$-categorical. Thus the
  countable models of $\Sigma$ are just the countable models of
  $\Cat$. Hence $\Mod(\Sigma)$ and each $\D_\kappa$ satisfy II and
  III. The infinite-dimensional model $M_\omega$ lies in $\Mod(\Sigma)$ and
  each $\D_\kappa$, so IV.2 holds. It remains to prove axioms I.2 and
  IV.1. We prove two related families of statements:
  \begin{tightlist}
  \item[1${}_\kappa$)] $\D_\kappa$ satisfies I.2 and hence is a
    \qmeclass.
  \item[2${}_\kappa$)] If $H \models \Sigma$ and $N$ is a closed
    subset of $H$ with $\dim N = \kappa$ then $N \models \Sigma$.
  \end{tightlist}
  We must prove statement 2 for all $\kappa$, and statement 1 for
  uncountable $\kappa$. First we prove statement 2 for countable
  $\kappa$. For $n\in\N$, let $x_1,\ldots,x_n$ be variables not
  occuring in $\sigma_n$, and let $\sigma_n'(x_1,\ldots,x_n)$ be the
  $\Loo$-formula obtained from $\sigma_n$ by recursively replacing all
  quantified subformulas of the form $\exists y [\phi(y)]$ by $\exists
  y [\pi_n(y,x_1\ldots,x_n) \wedge \phi(y)]$, for any subformula
  $\phi(y)$ (and similarly for universal quantifiers). Let
  $\Indep_n(x_1,\ldots,x_n)$ be the formula $\bigwedge_{i=1}^n
  \neg\pi_{i-1}(x_i,x_1,\ldots,x_{i-1})$. Then
  $\Indep_n(x_1,\ldots,x_n)$ says that the $x_i$ are
  $\cl$-independent, and $\sigma_n'(x_1,\ldots,x_n)$ says that the
  closure of the $x_i$ is a model of $\sigma_n$. Thus for each $m \le
  \omega$, and each $n \in \N$, 
  \[\sigma_m \entails (\forall x_1,\ldots,x_n)[\Indep_n(x_1,\ldots,x_n) \to
  \sigma_n'(x_1,\ldots,x_n)].\]
  Hence if $N$ is finite dimensional, it is a model of $\Sigma$.

  Now suppose $\dim N = \aleph_0$. Let $(b_n)_{n < \omega}$ be a basis
  for $N$, and let $N_m = \cl_H(\class{b_n}{n<m})$ for each $m <
  \omega$. Then $N$ is the union of the chain $(N_m)_{m < \omega}$
  and, by the above, $N_m \iso M_m$. The union of the chain
  $(M_m)_{m\le \omega}$ is $M_\omega$, hence $N \iso M_\omega$. In
  particular, $N \models \Sigma$. Thus $2_\kappa$ holds for all
  countable $\kappa$.

  Now we prove statements $1_\kappa$ and $2_\kappa$ together for
  uncountable $\kappa$, by induction. Suppose inductively that
  $2_\lambda$ holds for every $\lambda < \kappa$. Then if $H \in
  \D_\kappa$, every closed subset of $H$ has dimension less than
  $\kappa$, so is a model of $\Sigma$. Hence $1_\kappa$ holds.

  Now suppose $H \models \Sigma$ and $N$ is a closed subset of $H$
  with $\dim N = \kappa$. Identifying $\kappa$ with its initial
  ordinal, let $(b_\lambda)_{\lambda < \kappa}$ be a basis of $N$ and
  let $N_\mu = \cl_H(\class{b_\lambda}{\lambda < \mu})$. Then $N$ is
  the union of the chain $(N_\mu)_{\omega \le \mu < \kappa}$, and by
  induction each $N_\mu$ models $\Sigma$. The chain lies in
  $\D_\kappa$, which by the same induction hypothesis satisfies I.2
  and hence is a \qmeclass.  Thus by lemma~\ref{Loo-embedding}, the
  chain is an $\Loo$-chain.  Each $N_\mu$ in the chain is infinite
  dimensional, thus models $\sigma_\omega$, and hence $N \models
  \sigma_\omega$. Also $N \subs H$, and $H$ has the CCP, so $N$ also
  has the CCP. Thus $N \models \Sigma$, that is, $2_\kappa$ holds.
  Thus, by induction, $1_\kappa$ and $2_\kappa$ hold for all
  uncountable $\kappa$.

  Thus $\Mod(\Sigma)$ satisfies axiom I.2, and so is a quasiminimal
  excellent class. If $(H_\lambda)_{\lambda<\kappa}$ is any chain in
  $\Mod(\Sigma)$ then either its union is finite dimensional (and lies
  in $\Mod(\Sigma)$) or the chain is eventually infinite dimensional
  and by lemma~\ref{Loo-embedding} is eventually an $\Loo$-chain. In
  the latter case, as above, the union $H$ of the chain is a model of
  $\sigma_\omega$. If $X$ is a finite subset of $H$ then $X \subs
  H_\lambda$ for some $\lambda < \kappa$, and $\cl_H(X) =
  \cl_{H_\lambda}(X)$. Since $H_\lambda$ has the CCP, this closure is
  countable. Hence $H$ has the CCP, and so $H \models \Sigma$. Thus
  $\Mod(\Sigma)$ satisfies axiom IV.1. By theorem~\ref{existence},
  $\Mod(\Sigma)$ is an uncountably categorical \qmeclass.

  If $H \in \Cat$ then either $H$ is finite dimensional in which case
  $H$ is $M_n$ for some $n \in \N$, or there is a closed embedding
  $M_\omega \into H$. By lemma~\ref{Loo-embedding}, this embedding is
  an $\Loo$-embedding, so $H \models \sigma_\omega$. In either case,
  $H \models \Sigma$. So $\Cat \subs \Mod(\Sigma)$. Since $\Cat$
  satisfies I.2, it is an initial segment of $\Mod(\Sigma)$. That is,
  either $\Cat = \D_\kappa$ for some $\kappa$ or $\Cat$ satisfies IV.1
  and $\Cat = \Mod(\Sigma)$.
\end{proof}

\begin{cor}\label{unique C'}
  Let $\Cat$ be a \qmeclass\ with a model of dimension
  $\aleph_0$. Then there is a unique \qmeclass\ $\Cat'$, containing
  $\Cat$, which is uncountably categorical.
\end{cor}
\begin{proof}
  By proposition~\ref{countable language} we may assume the language
  is countable. Then theorem~\ref{Scott sentence} gives $\Mod(\Sigma)$
  as one such $\Cat'$. If $\Cat''$ is any uncountably categorical
  \qmeclass\ containing $\Cat$ then by theorem~\ref{existence} it
  satisfies axiom IV, so by theorem~\ref{Scott sentence} again it is
  equal to $\Mod(\Sigma)$. Hence $\Cat'$ is unique.
\end{proof}
From this result we see that nothing would be lost by adding axiom IV
to the definition of a \qmeclass, unless perhaps there can be
interesting behaviour of finite dimensional models. See the questions
at the end of this paper.

In theorem~\ref{Scott sentence}, the sentence $\Sigma$ depends only on
the model $M_\omega$, because the $M_n$ are substructures of
$M_\omega$. We can extract the properties of $M_\omega$ which are
needed to produce a \qmeclass.
\begin{cor}\label{bootstrap}
  Let $M$ be a countable $\L$-structure, equipped with a pregeometry
  $\cl$, satisfying the following axioms.
  \begin{description}
    \item[I$'$ (Pregeometry)] \ \\ 
      The pregeometry is quantifier-free $\Loo$-definable, and $\dim M
      = \aleph_0$.

    \item[II$'$ ($\aleph_0$-homogeneity over closed sets)] \ \\ Let $G
      \subs M$ be closed or empty.
      \begin{tightlist}
        \item[]
        \item[II.1$'$] If $x,x' \in M$ are each independent from $G$,
          then $\qftp(x/G) = \qftp(x'/G)$.
	\item[II.2$'$] Let $\bar{x},\bar{x}'$ be finite tuples from
          $M$ such that $\qftp(\bar{x}/G) = \qftp(\bar{x}'/G)$, and
          let $y \in M$. Then there is $y'\in M$ such that
          $\qftp(\bar{x}y/G) = \qftp(\bar{x}'y'/G)$.
      \end{tightlist}

    \item[III$'$ (Excellence)]\ \\ 
      If $C$ is a crown in $M$ and $\bar{x}$ is a finite tuple from
      $\cl(C)$, then there is a finite subset $C_0$ of $C$ such that
      for any tuple $\bar{x}'$, 
      \[\qftp(\bar{x}/C_0) = \qftp(\bar{x}'/C_0) \implies
      \qftp(\bar{x}/C) = \qftp(\bar{x}'/C). \]
  \end{description}
  Then there is a unique \qmeclass\ $\Cat$ which satisfies axiom IV
  such that $\tuple{M,\cl} \in \Cat$.
\end{cor}

\begin{proof}
  Let $\Cat_0$ be the class of $\L$-structures consisting of $M$ and
  all its $cl$-closed substructures, equipped with the pregometry
  $\cl$ and its restrictions. The axioms I$'$---III$'$ ensure that
  $\Cat_0$ is a \qmeclass\ with a model of dimension $\aleph_0$, and
  then corollary~\ref{unique C'} gives the unique class $\Cat$.
\end{proof}

\section{Questions}
We conclude with some further questions.
\begin{enumerate}
  \item Is there a class which satisfies axioms I and II, but not III?
    It seems likely that there is, but I do not believe that any
    examples are currently known.

  \item Does axiom III follow from lemma~\ref{excellence lemma}? The
    former is Zilber's definition of excellence and the latter should
    be Shelah's with respect to this abstract elementary class, albeit
    with prime model in place of primary model.

  \item If $\Cat$ is a \qmeclass\ with models of arbitrarily large
    finite dimension (and hence of all finite dimensions), there is a
    well-defined $\L$-structure $M$ constructed as the union of all of
    the finite dimensional models. Is $\Cat \cup \{M\}$ necessarily a
    \qmeclass?

  \item What \qmeclass es are there with models only of dimension up
    to some finite $n$, with $n \ge 1$? 
\end{enumerate}

\bibliographystyle{../pjk1}
\bibliography{../papers,../books,../drafts}

\end{document}